\documentclass[11pt,twoside]{article}
\title{} \author{} \date{}
\usepackage{latexsym,amssymb,times}
\input amssym.def
\newtheorem{te}{Theorem}[section]
\newtheorem{prop}[te]{Proposition}

\newtheorem{fac}[te]{Fact}
\newtheorem{lem}[te]{Lemma}

\newtheorem{ex}[te]{Example}

\newtheorem{cl}[te]{Claim}

\def\dok{\noindent{\bf Proof. }}
\def\kdok{\hfill $\Box$ \par \vspace*{2mm} }
\def\a{\alpha}
\def\b{\beta}
\def\k{\kappa}
\def\l{\lambda}

\def\o{\omega}
\def\t{\tau}
\def\s{\sigma}
\def\f{\varphi}
\def\p{\psi}
\def\fc{{\mathfrak c}}
\def\N{{\mathbb N}}
\def\Q{{\mathbb Q}}
\def\R{{\mathbb R}}
\def\Z{{\mathbb Z}}
\def\I{{\mathbb I}}

\def\BL{{\mathbb L}}
\def\BA{{\mathbb A}}
\def\BC{{\mathbb C}}
\def\BD{{\mathbb D}}

\def\O{{\mathcal O}}

\def\L{{\mathcal L}}

\def\CS{{\mathcal S}}
\def\CT{{\mathcal T}}

\def\la{\langle}
\def\ra{\rangle}

\def\Top{\mathop{{\mathrm{Top}}}\nolimits}

\def\cf{\mathop{{\mathrm{cf}}}\nolimits}
\def\ci{\mathop{{\mathrm{ci}}}\nolimits}

\def\id{\mathop{{\mathrm{id}}}\nolimits}
\def\Aut{\mathop{{\mathrm{Aut}}}\nolimits}

\def\Sym{\mathop{{\mathrm{Sym}}}\nolimits}
\def\Conv{\mathop{{\mathrm{Conv}}}\nolimits}

\begin{document}
\thispagestyle{plain}
\begin{center}
           {\large \bf \uppercase{Partial orders of bijectively related or\\[1mm]
                        homeomorphic topologies}}
\end{center}
\begin{center}
{\bf Aleksandar Janjo\v s\footnote{Faculty of Technical Sciences, University of Novi Sad,
                                      Trg Dositeja Obradovi\'ca 6, 21000 Novi Sad, Serbia.
                                      e-mail: aleksandar.janjos@uns.ac.rs}
and Milo\v s S.\ Kurili\'c}\footnote{Department of Mathematics and Informatics, Faculty of Sciences, University of Novi Sad,
                                      Trg Dositeja Obradovi\'ca 4, 21000 Novi Sad, Serbia.
                                      e-mail: milos@dmi.uns.ac.rs}
\end{center}
\begin{abstract}
\noindent
Topologies $\tau , \sigma \in \mathop{{\mathrm{Top}}}\nolimits _X$
are {\it bijectively related}, in notation $\tau \sim \sigma$,
if there are continuous bijections $f: (X, \tau )\rightarrow  (X, \sigma )$ and $g: (X, \sigma)\rightarrow  (X, \tau)$.
Defining $[\tau ]_{\cong}=\{ \sigma \in \mathop{{\mathrm{Top}}}\nolimits _X : \sigma \cong \tau\}$ 
and $[\tau ]_{\sim }=\{ \sigma \in \mathop{{\mathrm{Top}}}\nolimits _X : \sigma \sim \tau\}$
we show that for each infinite 1-homogeneous linear order ${\mathbb L}$
there is a topology $\tau \in \mathop{{\mathrm{Top}}}\nolimits _{|L|}$ such that
\begin{itemize}
\item[(a)] $\langle [\tau ]_{\cong}, \subset \rangle \cong \dot{\bigcup}_{2^{|L|}}{\mathbb L}$ (the disjoint union of $2^{|L|}$-many copies of ${\mathbb L}$);\\
so, each maximal chain in $[\tau ]_{\cong}$ is isomorphic to ${\mathbb L}$;
\item[(b)] $\langle [\tau ]_{\sim}, \subset \rangle\cong \dot{\bigcup}_{2^{|L|}}\widetilde{{\mathbb L}}$, where $\widetilde{{\mathbb L}}$ is the Dedekind completion of ${\mathbb L}$;\\
thus, each maximal chain in $[\tau ]_{\sim}$ is isomorphic to $\widetilde{{\mathbb L}}$.
\end{itemize}
If, in addition, the linear order ${\mathbb L}$ is Dedekind complete, then the topology $\tau$ is weakly reversible, non-reversible
and $\langle [\tau ]_{\sim}, \subset \rangle=\langle [\tau ]_{\cong}, \subset \rangle\cong \dot{\bigcup}_{2^{|L|}}{\mathbb L}$.



\end{abstract}
\section{Introduction}\label{S1}
The partial order $\la \Top _X ,\subset\ra$, where $\Top _X$ is the set of all topologies on a non-empty set $X$,
was widely considered; it is a complete lattice:
$\bigwedge \CT =\bigcap \CT$ and $\bigvee \CT =\t(\bigcup \CT)$, for a family of topologies $\CT \subset \Top _X$,
where $\t (\CS)$ denotes the minimal topology on the set $X$ containing a family $\CS$ of its subsets.

In fact, regarding the topological structures on a set $X$,
we are interested for the diversity of {\it non-homeomorphic} topologies on $X$;
i.e., for the structure of the quotient $\Top _X \!/\!\cong$,
while, assuming in the sequel that the set $X$ is infinite,
each topology $\t \in \Top _X$ (except a few trivial ones) has infinitely many incarnations in $\Top _X$
(the copies of $\t$, $f[\t]:=\{ f[O]:O\in \t\}$, obtained by the action of the symmetric group $\Sym (X )$).
Therefore, it is reasonable to compare the equivalence classes $[\t ]_{\cong}:=\{ \s \in \Top _X : \s \cong \t\}$
(or topologies $\t\in \Top _X$)
saying that $[\t _1]_{\cong} \preccurlyeq [\t _2]_{\cong}$  (or $\t _1 \preccurlyeq \t _2)$ iff there is $\t \in \Top _X$ such that $\t _1 \cong \t \subset \t _2$,
or, equivalently, iff there is a continuous bijection (a condensation) $f: \t _2\rightarrow \t _1$ (more precisely, $f: (X, \t _2)\rightarrow  (X, \t _1)$).
So we obtain a preorder $\la \Top _X \!/\!\cong ,\preccurlyeq \ra $
and, defining the {\it condensational equivalence} relation by $\t _1 \sim \t _2$ iff $\t _1 \preccurlyeq \t _2$ and $\t _2 \preccurlyeq \t _1$,
(such topologies are called {\it bijectively related} in \cite{Doy})
we obtain its antisymmetric quotient, the {\it condensational order} $\la \Top _X \!/\!\sim ,\leq \ra$,
where $[\t _1]_{\sim} \leq [\t _2]_{\sim}$ iff $\t _1 \preccurlyeq \t _2$.

For each topology $\t\in \Top _X$ we have $[\t ]_{\cong}\subset  [\t ]_{\sim}=\Conv([\t ]_{\cong})$,
where $\Conv (\CT)$ denotes the convex closure of a suborder $\CT$ of the lattice $\la \Top _X ,\subset\ra$,
and we distinguish the following four cases (where 1 implies 2 implies 3; see \cite{K}).

{\it Case 1}: $[\t ]_{\sim}=[\t ]_{\cong}=\{ \t\}$.
We call such topologies {\it strongly reversible};
here we only have the discrete topology, the antidiscrete topology and generalizations of the cofinite topology.

{\it Case 2}: $[\t ]_{\sim}=[\t ]_{\cong}$ is an antichain in the lattice $\la \Top _X ,\subset\ra$.
These are exactly the {\it reversible} topologies: each continuous bijection $f: \t \rightarrow \t $ is a homeomorphism; see \cite{Raj}.
We recall that the class of reversible spaces includes several prominent structures as the Euclidean spaces $\R ^n$ (by Brouwer's invariance of domain theorem)
and compact Hausdorff spaces, while many prominent structures are not reversible e.g.\ $\Q$, $\R \setminus \Q$ and normed spaces of infinite dimension (see \cite{Raj}).

{\it Case 3}: $[\t ]_{\sim}=[\t ]_{\cong}$ is a convex set in the lattice $\la \Top _X ,\subset\ra$. These are {\it weakly reversible} topologies in our terminology.
These are exactly the topologies satisfying the Cantor-Schr\"{o}eder-Bernstein property for condensations:
if there are condensations $f:\s \rightarrow \t$ and $g:\t \rightarrow \s$, there is a homeomorphism  $F:\s \rightarrow \t$.

{\it Case 4}: $[\t ]_{\sim}\supsetneq [\t ]_{\cong}$.
Such topologies are not weakly reversible
and here we have $[\t ]_{\sim}=\dot{\bigcup}_{\s \sim \t} [\s ]_{\cong}$,
where the condensation class $[\t ]_{\sim}$ is a union of at least two homeomorphism classes.

If instead of the complete lattice $\la \Top _X ,\subset\ra$ we consider the partial order $\la \Top _X \!/\!\sim ,\leq \ra$
(which is not a lattice\footnote{e.g., if $x,y,z\in X$ are different points,
then $[\{ \emptyset ,\{ x\}, X\}]_{\sim} \lor [\{ \emptyset ,\{ x,y\}, X\}]_{\sim}$ does not exist}),
its suborders and the equivalence classes $[\t ]_{\sim}$ and $[\t ]_{\cong}$,
then several questions arise.
For example, Rajagopalan and Wilansky asked
what is the cardinality of a maximal chain of homeomorphic topologies (see \cite{Raj}, p.\ 137).
Some natural extensions of that question are

- Which are the possible order types of maximal chains in 
$[\t ]_{\cong}$ and $[\t ]_{\sim}$?

- How do the partial orders $\la [\t ]_{\cong}, \subset \ra$ and $\la [\t ]_{\sim}, \subset \ra$ look like?

\noindent
In this paper we consider non-reversible topologies.
It is easy to see that, for such a topology $\t$,
maximal chains in $[\t ]_{\cong}$ have no end points
and maximal chains in $[\t ]_{\sim}$ are Dedekind complete.
In addition, a weakly reversible topology of a sequential space in which each sequence has at most one limit must be reversible \cite{K};
thus, roughly speaking, ``nice" (e.g.\ metrizable) weakly reversible, non-reversible topologies do not exist.
In Section \ref{S3} we show that for each 1-homogeneous linear order $\BL $ of size $\k\geq \o$
there is a topology $\t \in \Top _\k$ such that
$\la [\t ]_{\cong}, \subset \ra\cong \dot{\bigcup}_{2^\k}\BL$,
(where, generally, $\dot{\bigcup}_{\l}\BL$ denotes the disjoint union of $\l$-many copies of a linear order $\BL$)
so, each maximal chain in $[\t ]_{\cong}$ is isomorphic to $\BL$;
also,  $\la [\t ]_{\sim}, \subset \ra\cong \dot{\bigcup}_{2^\k}\widetilde{\BL}$, where $\widetilde{\BL}$ is the Dedekind completion of $\BL$;
thus, each maximal chain in $[\t ]_{\sim}$ is isomorphic to $\widetilde{\BL}$.
If the linear order $\BL$ is, in addition, Dedekind complete, then the topology $\t$ is weakly reversible
and $\la [\t ]_{\sim}, \subset \ra=\la [\t ]_{\cong}, \subset \ra\cong \dot{\bigcup}_{2^\k}\BL$.
Otherwise, for $\s\in [\t ]_{\sim}\setminus [\t ]_{\cong}$ there is a linear order $\BL _\s$ such that $\la [\s ]_{\cong}, \subset \ra\cong \dot{\bigcup}_{2^\k}\BL_\s$.
\section{Preliminaries}\label{S2}
We recall that a linear order $\BL =\la L,<\ra$ is {\it Dedekind complete}
iff each non-empty set $X\subset L$ having an upper bound has a least upper bound, $\sup X $.
A {\it Dedekind completion} of a linear order $\BL $ is the unique (up to isomorphism) minimal Dedekind complete extension $\widetilde{\BL}=\la\widetilde{L}, \prec\ra$ of $\BL$
and one of its isomorphic incarnations is obtained by adding to $\BL$ the set of its gaps and extending $<$ to $\prec$ in a natural way;
a {\it gap} in $\BL$ is a pair $\la I,F\ra$, where $\{ I,F\}$ is a partition of $L$ into non-empty sets, $I<F$,
$I$ has no last element and  $F$ has no first element.
The minimality of $\widetilde{\BL}$ can be expressed in the following way: for each $c,c'\in \widetilde{L}$

(m1) $L\cap (\cdot ,c]\neq \emptyset\;$ and  $\;L\cap [c, \cdot )\neq \emptyset$,

(m2) If $c \prec c'$, then $| L\cap[c,c']|\geq 2$.

\noindent
(E.g.\ for (m2): if $c,c'\in L$, we are done; if $c\not\in L$, then $c=\la I,F\ra$ and $|(c,c')\cap F|\geq \o$, etc.) The following fact will be used in the sequel.
\begin{fac}\label{T027}
If $\BL$ is a linear order, $\widetilde{\BL}$ its Dedekind completion and $\BL\subset\widetilde{\BL}$, then

(a) If $\BL$ has no end-points, then $\widetilde{\BL}$ has no end-points and $L$ is a cofinal and coinitial set in $\widetilde{\BL}$;

(b) If  $c,c'\in \widetilde{L}$ and $c\prec c'$, then $(\cdot,c)\cap L\subsetneq(\cdot,c')\cap L$.

(c) Each automorphism of $\BL$ has a unique extension to an automorphism of $\widetilde{\BL}$.

\end{fac}
\dok
(a) If $c\in \widetilde{L}$, then by (m1) there is $x\in L\cap [c, \cdot )$
and, since $\BL$ has no last element,  there is $x'\in L$ such that $x'>x \geq c$.
So the set $L$ is cofinal in $\widetilde{\BL}$ and $\widetilde{\BL}$ has no last element. The rest of the proof is dual.

(b) By (m2)  there is $x\in L\cap [c,c')$ and, hence, $x\in(\cdot,c')\cap L$ and $x\not\in(\cdot,c)\cap L$.

(c) If $f\in\Aut (\BL)$, it is easy to check that $f\subset F\in \Aut (\widetilde{\BL})$, where $F(c)=\sup _{\widetilde{\BL}} f[L \cap (\cdot ,c]]$, for $c\in\widetilde{L}$.
\kdok
A linear order $\BL$ is called {\it 1-homogeneous} (or 1-transitive) iff for each $x,y\in L$
there is an automorphism $f \in \Aut (\BL)$ such that $f(x)=y$.
Then, clearly, $|L|>1$ implies that $L$ has no end points.
More generally, $\BL$ is called {\it $n$-homogeneous} (or $n$-transitive) iff for each $k\leq n$ and each $x_1,\dots ,x_k,$ $y_1,\dots, y_k\in L$,
such that $x_1<x_2<\dots <x_k$ and $y_1<y_2<\dots <y_k$,
there is an automorphism $f \in \Aut (\BL)$ such that $f(x_i)=y_i$, for all $i\leq k$.
$\BL$ is called {\it ultrahomogeneous} if it is $n$-homogeneous, for each $n\in \N$. We recall some basic facts (see \cite{Ros}).
\begin{fac}\label{T012}
If $\BL $ is a linear order, then we have

(a) $\BL $ is 1-homogeneous iff $\BL \cong \Z ^\a \BD$, where $\a$ is an ordinal and $\BD$ is an 1-homogeneous dense linear order;

(b) $\BL $ is 1-homogeneous, scattered iff $\BL \cong \Z ^{r_F(\BL)}$;

(c) $\BL $ is 1-homogeneous, countable iff $\BL \cong \Z ^\a$ or $\BL \cong \Z ^\a \Q$, where $\a<\o _1$;

(d) $\BL $ is 2-homogeneous  iff it is $n$-homogeneous for all $n\geq 2$;

(e) If $\BL $ and $\BL '$ are 1-homogeneous linear orders, then $\BL \BL '$ is 1-homogeneous.
\end{fac}
\begin{ex}\label{EX000}\rm
It is evident that the linear orders $\Z$, $\Q$, $\R$ and $\I =\R\setminus \Q$ are 1-homogeneous.
$\Z$ is not 2-homogeneous and $\Z +\Q$ is not 1-homogeneous.
Taking $\Z$, $\Q$, $\R$ and $\I$, and their $4\times 4=16$ lexicographical products, $\Z\Z ,\Z\Q,\dots, \I\I$,
by Fact \ref{T012}(e) we obtain 1-homogeneous linear orders.
\end{ex}
\section{Non-reversible topologies}\label{S3}
\begin{te}\label{T003}
For each $1$-homogeneous linear order $\BL $ of size $\k\geq \o$ there is a non-reversible topology $\t \in \Top _\k$ such that

(a) $\la [\t ]_{\cong}, \subset \ra\cong \dot{\bigcup}_{2^\k}\BL$; so, each maximal chain in $[\t ]_{\cong}$ is isomorphic to $\BL$;

(b) $\la [\t ]_{\sim}, \subset \ra\cong \dot{\bigcup}_{2^\k}\widetilde{\BL}$; so, each maximal chain in $[\t ]_{\sim}$ is isomorphic to $\widetilde{\BL}$;

(c) For $\s\in [\t ]_{\sim}\setminus [\t ]_{\cong}$ there is a linear order $\BL _\s$ such that $\la [\s ]_{\cong}, \subset \ra\cong \dot{\bigcup}_{2^\k}\BL_\s$.
\end{te}
The rest of this section is devoted to a proof of the theorem.
Let $\BL =\la L,< _\BL\ra$ be a 1-homogeneous linear order of size $\k>1$,
$\widetilde{\BL} =\la \widetilde{L}, <_{\widetilde{\BL}}\ra$ its Dedekind completion,
where $\BL$ is a suborder of $\widetilde{\BL}$, and let $z\not\in \widetilde{L}$.
We regard the linear order
\begin{equation}\label{EQ003}
\widetilde{\BL} _z =\la \widetilde{L}_z , <\ra, \;\;\mbox{ where }\;\;\widetilde{L}_z:=\{ z\}\cup \widetilde{L} \;\;\mbox{ and }\;\; < \;:= \{ \la z,a\ra :a\in \widetilde{L}\} \;\cup <_{\widetilde{\BL}}
\end{equation}
(so, $z=\min (\widetilde{\BL}_z)$ and $\widetilde{\BL}_z=\{ z\}+\widetilde{\BL}$).
Clearly, $\BL _z :=\la L_z ,<\upharpoonright L_z  \ra$, where $L_z:= \{ z\}\cup L$,
is a suborder of the order $\widetilde{\BL}_z$.
For simplicity, by $(s,t)$, $[s,t)$, etc., where $s,t\in \widetilde{L}_z$,
we denote the intervals {\it in the linear order $\widetilde{\BL}_z$}; for example,
$[s,t)=\{ u\in \widetilde{L}_z : s\leq u<t\}$.
For $c\in \widetilde{L}$ let $\t _c :=  \t _\leftarrow \cup \s _c \cup \{ \emptyset , L_z\}$, where
\begin{equation}\label{EQ004}
\t _\leftarrow := \{ [z,a) \cap L_z : a\in \widetilde{L}\} \;\;\mbox{and}\;\; \s _c :=  \{ (z,b) \cap L_z :b\in (z,c]\},
\end{equation}
and let $x_0$ be an element of $L$.
We will show that the topology $\t _{x_0}\in \Top _{L_z}$ is as stated in Theorem \ref{T003} (see Propositions \ref{T037},  \ref{T041} and  \ref{T048}).
\subsection{A maximal chain of homeomorphic topologies}
\begin{prop}\label{T042}
$\L :=\{ \t _x :x\in L\}$ is a maximal chain in the poset $\la [\t _{x_0}]_{\cong}, \subset\ra$
isomorphic to $\BL$.
\end{prop}
\dok
It is easy to check that $\t_c$, $c\in\widetilde{L}$, are topologies on the set $L_z$.
E.g.\ if $A\subset \widetilde{L}$ is a non-empty set having an upper bound,
then by the Dedekind completeness of $\widetilde{\BL}$ we have $\bigcup_{a\in A}[z,a) \cap L_z =[z,a') \cap L_z \in \t _c$,
where $a':=\sup_{\widetilde{\BL}} A$, etc.
\begin{cl}\label{T016}
For each $c_1,c_2\in\widetilde{L}$ we have: $c_1<c_2\Leftrightarrow \t_{c_1}\subsetneq\t_{c_2}$.
\end{cl}
\dok
If $c_1<c_2$, then by Fact \ref{T027}(b) there is $x\in (\cdot,c_2)\setminus (\cdot,c_1) \cap L$,
which gives  $x\in (z,c_2)\cap L_z \in \s _{c_2}$.
Since $x\geq c_1$ we have $(z,c_2)\cap L_z \not\in \s _{c_1}$
and, hence, $\s_{c_1}\subsetneq\s_{c_2}$ and $\t_{c_1}\subsetneq\t_{c_2}$. The converse is evident.
\hfill $\Box$
\begin{cl}\label{T018}
(a) $\t_{x_1}\cong \t_{x_2}$, for all $x_1,x_2\in L$.

(b) $\L =\{ \t _x :x\in L\}$ is a chain in $\la [\t _{x_0}]_{\cong}, \subset\ra$ and $\la\L,\subsetneq\ra\cong \BL$.

(c) $\Lambda =\{ \t _c :c\in \widetilde{L}\}$ is a chain in $\la \Top _{L_z}, \subset\ra$,
$\la\Lambda,\subsetneq\ra\cong\widetilde{\BL}$ and $\Lambda =\widetilde{\L}$.
\end{cl}
\dok
(a) Since $\BL$ is 1-homogeneous there is $\varphi\in\Aut(\BL)$ such that $\varphi(x_1)=x_2$,
and by Fact \ref{T027}(c) there is $\Phi\in\Aut(\widetilde{\BL})$, where $\Phi|L=\varphi$;
so, $F:= \{\la z,z\ra\} \cup \Phi\in\Aut(\widetilde{\BL}_z )$.
Clearly $f:=F|L_z= \{ \la z,z\ra\} \cup \varphi \in\Aut(\BL_z)$
and we show that $f:(L_z,\t_{x_1})\rightarrow(L_z,\t_{x_2})$ is a homeomorphism.

If $a\in\widetilde{L}$, then since $F[L_z]=L_z$ we have
$f[[z,a)\cap L_z]=F[[z,a)\cap L_z]=[z,F(a))\cap L_z$; so
\begin{equation}\label{EQ000}
f[[z,a)\cap L_z]=[z,F(a))\cap L_z.
\end{equation}
For $b\in(z,x_1]$ we have $F(b)\in (F(z),F(x_1)]=(z,x_2]$ and
$f[(z,b)\cap L_z]= F[(z,b)\cap L_z]=(z,F(b))\cap L_z$; thus
\begin{equation}\label{EQ001}
f[(z,b)\cap L_z]=(z,F(b))\cap L_z \;\;\mbox{and}\;\; F(b)\in(z,x_2].
\end{equation}
By (\ref{EQ000}) and (\ref{EQ001}) the mapping $f$ is open and
we show that it is continuous.
For $a\in\widetilde{L}$ by $(\ref{EQ000})$ we have $f[[z,F^{-1}(a))\cap L_z]=[z,a)\cap L_z$,
and, hence, $f^{-1}[[z,a)\cap L_z]=[z,F^{-1}(a))\cap L_z \in \t _\leftarrow$.
For $b\in(z,x_2]$ we have  $F^{-1}(b)\in (F^{-1}(z),F^{-1}(x_2)]=(z,x_1]$
and by $(\ref{EQ001})$ we have $f[(z,F^{-1}(b))\cap L_z]=(z,b)\cap L_z$,
which implies that $f^{-1}[(z,b)\cap L_z]=(z,F^{-1}(b))\cap L_z \in \tau _{x_1}$.
Thus $f$ is continuous indeed.

(b) By (a) we have $\L \subset [\t _{x_0}]_{\cong}$.
By Claim \ref{T016} the mapping $f:\BL\rightarrow \la\L,\subsetneq\ra$ defined by $f(x)=\t_x$, for $x\in L$, is an isomorphism.

(c) By Claim \ref{T016} the mapping $F:\widetilde{\BL}\rightarrow\la\Lambda ,\subsetneq\ra$ defined by $F(c)=\t_c$, for $c\in \widetilde{L}$, is an isomorphism.
Thus the linear order $\Lambda$ is Dedekind complete and, since $F\upharpoonright \L =f$, $\L$ is dense in $\Lambda$; so $\Lambda$ is a Dedekind completion of $\L$.
\hfill $\Box$
\begin{cl}\label{T020}
For $c\in \widetilde{L}$, $z$ is the unique point $x\in L_z$ such that $\bigcap _{ x\in O\in\t_c} O=\{x\}$.
\end{cl}
\dok
Since $\widetilde{\BL}$ has no first element we have $\bigcap _{ z\in O\in\t_c} O=\bigcap_{a\in\widetilde{L}}[z,a)\cap L_z=\{z\}$.
If $x>z$, then, since $\BL$ has no first element, there is $y\in (z,x)\cap L_z$
and $x\in O\in\t_c$ implies $y\in O$;
thus $y\in \bigcap _{x\in O\in\t_c } O\neq\{x\}$.
\hfill $\Box$
\begin{cl}\label{T021}
If $c_1,c_2\in \widetilde{L}$ and $f:(L_z,\t_{c_1})\rightarrow (L_z,\t_{c_2})$ is a homeomorphism, then

(a) $f(z)=z$;

(b) $f[(z,c_1)\cap L_z]=(z,c_2)\cap L_z$.
\end{cl}
\dok
(a) For $O\subset L_x$ we have $z\in O\in\t_{c_1}$ iff $f(z)\in f[O]\in\t_{c_2}$;
so, by Claim \ref{T020}, $f(\{z\})=f[\bigcap_{z\in O\in\t_{c_1}} O]=\bigcap_{f(z)\in f[O]\in\t_{c_2}}\ f[O] \supset\bigcap_{f(z)\in O'\in\t_{c_2}} O'\ni f(z)$.
Thus $\bigcap _{f(z)\in O'\in\t_{c_2}}O'=\{ f(z)\}$
and, by Claim \ref{T020} again, $f(z)=z$.

(b) Since $z\notin (z,c_1)\cap L_z$
by (a) we have $z=f(z)\not\in f[(z,c_1)\cap L_z]$
and, hence, $f[(z,c_1)\cap L_z] \in\s_{c_2}$,
which gives $f[(z,c_1)\cap L_z]\subset (z,c_2)\cap L_z$.
Since $f^{-1}$ is a homeomorphism too,
in a similar way we obtain $f^{-1}[(z,c_2)\cap L_z]\subset (z,c_1)\cap L_z$
and, thus, $(z,c_2)\cap L_z\subset f[(z,c_1)\cap L_z]$.
\hfill $\Box$
\begin{cl}\label{T022}
If $c\in\widetilde{L}\setminus L$, then $\t_{c}\not\cong\t_{x_0}$.
\end{cl}
\dok
Suppose $f:(L_z,\t_{c})\rightarrow(L_z,\t_{x_0})$ is a homeomorphism.
Since $(z,c]\cap L_z=(z,c)\cap L_z$
and, by Claim \ref{T021}(b), $f[(z,c)\cap L_z]=(z,x_0)\cap L_z$
by Claim \ref{T021}(a) we have
\begin{equation}\label{EQ002}
f[[z,c]\cap L_z]=[z,x_0)\cap L_z.
\end{equation}
Let $x\in L$, where $f(x)=x_0$.
Since $c\not\in L$ by (\ref{EQ002}) we have $x>c$.
By (m2) we have $|L\cap[c,x]|\geq 2$
and, since $c\not\in L$, there is $y\in(c,x)\cap L_z$.
Since $c\not\in L$, by Fact \ref{T027}(b) we have $[z,c]\cap L_z =[z,c)\cap L_z\subsetneq[z,y)\cap L_z$
which by (\ref{EQ002}) implies that
$[z,x_0)\cap L_z=f[[z,c]\cap L_z]\subsetneq f[[z,y)\cap L_z]=[z,a)\cap L_z$, for some $a\in\widetilde{L}$,
because $z\in[z,y)\cap L_z$ and by Claim \ref{T021}(a), $z\in f[[z,y)\cap L_z]$.
Consequently we have $x_0<a$;
thus $f(x)=x_0<a$
and, hence, $f(x)\in[z,a)\cap L_z$,
which implies that $x\in f^{-1}[[z,a)\cap L_z]=[z,y)\cap L_z$.
Since $x>y$ we have a contradiction.
\hfill $\Box$
\begin{cl}\label{T023}
If $A\subset\widetilde{L}$ is a non-empty set, then

(a) If $a=\sup _{\widetilde{\BL}} A$, then $\O(\bigcup_{c\in A}\t_c)=\t_a$;

(b) If  $a=\inf _{\widetilde{\BL}} A$, then $\bigcap_{c\in A}\t_c=\t_a$;

(c) $\bigcap_{x\in L}\t_x =\t_{\leftarrow}\cup\{\emptyset,L_z\}$;

(d) $\bigcup_{x\in L}\t_x=\t_{\leftarrow}\cup\{\emptyset,L_z\}\cup\{(z,b)\cap L_z:\ b\in\widetilde{L}\}$.
\end{cl}
\dok
(a)
The inclusion ``$\subset$" follows from Claim \ref{T016}.
Assuming that $\bigcup_{c\in A}\t_c\subset\t\in\Top_{L_z}$ we prove $\t_a\subset\t$.
Since $\t_{\leftarrow}\cup\{\emptyset,L_z\}\subset\t$, it remains to be shown that $\s_a=\{ (z,b)\cap L_z: b\in (z,a]\}\subset\t$.
Let $b\in(z,a]$.
If $b<a =\sup A$, then  $b<c$ for some $c\in A$,
and, hence, $(z,b)\cap L_z\in\t_c\subset\t$.
If $b=a$ and $x\in (z,a)\cap L_z$,
then $c_x>x$, for some $c_x\in A$
and, hence, $(z,a)\cap L_z \subset \bigcup _{x\in (z,a)\cap L_z}(z,c_x)\cap L_z \subset (z,a)\cap L_z$.
So, $(z,a)\cap L_z = \bigcup _{x\in (z,a)\cap L_z}(z,c_x)\cap L_z$
is a union of sets from $\bigcup_{c\in A}\t_c\subset\t$
and, since $\t$ is a topology, $(z,a)\cap L_z\in\t$.

(b)
The inclusion ``$\supset$" follows from Claim \ref{T016}.
Since $\bigcap_{c\in A}\t_c =\t_{\leftarrow}\cup\{\emptyset,L_z\}\cup\bigcap_{c\in A}\s_c$
it remains to be proved that $\bigcap_{c\in A}\s_c\subset\s_a$.
If $O\in\bigcap_{c\in A}\s_c$,
then for each $c\in A$ there is $b_c\in(z,c]$ such that $O=(z,b_c)\cap L_z$.
By Fact \ref{T027}(b) there is $b$ such that $b=b_c$, for all $c\in A$.
Thus $b\leq c$, for all $c\in A$,
and, since $a=\inf A$ we have $b\leq a$;
so, $O=(z,b)\cap L_z\in\s_a$.

(c) The inclusion ``$\supset$" is evident and for ``$\subset$" it remains to be proved that $\bigcap_{x\in L}\s_x =\emptyset$.
Assuming that $(z,b)\cap L_z \in \s_x$, for all $x\in L$, we would have $b\leq L$,
which is false since $\BL$ has no end points and the same holds for $\widetilde{\BL}$.

(d) The inclusion ``$\subset$" is evident.
If $b\in\widetilde{L}$, then, since $\BL$ has no end points, by Fact \ref{T027}(a) there is $x\in L$ such that $b\leq x$ and, hence, $(z,b)\cap L_z \in \t _x$.
\hfill $\Box$
\begin{cl}\label{T025}
$\L =\{ \t _x :x\in L\}$ is a maximal chain in the poset $\la [\t_{x_0}]_{\cong},\subset \ra$.
\end{cl}
\dok
Assuming that $\t\in[\t_{x_0}]_{\cong}$ and $\t\subsetneq\t_x$, for all $x\in L$,
by Claim \ref{T023}(c) we would have $\t \subset \t_{\leftarrow}\cup\{\emptyset,L_z\}$.
But then $\bigcap (\t \setminus \{ \emptyset\})=\{ z\}$, while $\bigcap (\t _{x_0}\setminus \{ \emptyset\})=\emptyset$;
thus, $\t \not\cong \t _{x_0}$ and we have a contradiction.

Assuming that $\t\in[\t_{x_0}]_{\cong}$ and $\t _x\subsetneq\t$, for all $x\in L$,
by Claim \ref{T023}(d) we would have $\t_{\leftarrow}\cup\{\emptyset,L_z\}\cup\{(z,b)\cap L_z:\ b\in\widetilde{L}\}\subset \t$
and, hence, $L=\bigcup _{b\in\widetilde{L}}(z,b)\cap L_z\in\t$.
Thus the singleton $\{ z\}$ would be a closed set in the space $\la L_z,\t\ra$,
while the space $\la L_z,\t _{x_0}\ra$ has no closed singletons.
So, $\t \not\cong \t _{x_0}$ and we have a contradiction.

Finally suppose that $\t\in[\t_{x_0}]_{\cong}\setminus\L$, where $\L\cup\{\t\}$ is a chain,
$I=\{x\in L:\ \t_x\subsetneq\t\}\neq\emptyset$ and $F=\{x\in L:\ \t\subsetneq\t_x\}\neq\emptyset$.
By Claims \ref{T016} and \ref{T018}(b) we have  $L=I+F$. Namely, $(I,F)$ is a cut in $\BL$ and we consider the following cases.

{\it Case 1}: $\max _{\BL}I=x'$ and $\min _{\BL}F=x''$.
Then $x'< x''$, $(x',x'')=\emptyset$ ($x'$ and $x''$ form a jump in $\BL$)
and, hence, $\t_{x''}\setminus\t_{x'}=\s_{x''}\setminus\s_{x'}=\{(z,x'')\cap L_z\}$,
which contradicts our assumption that $\t_{x'}\subsetneq\t\subsetneq\t_{x''}$.

{\it Case 2}: $\max _{\BL} I=x'$ and $\neg \exists \min _{\BL}F$.
Then $\t_{x'}\subsetneq\t$
and, since $x'=\inf F$, by Claim \ref{T023}(b) we have $\t_{x'}=\bigcap_{x\in F}\t_x\supset\t$,
which gives a contradiction.

{\it Case 3}: $\neg \exists\max _{\BL} I$ and $\min _{\BL}F=x''$.
Then $\t\subsetneq\t_{x''}$
and, since $x''=\sup I$, by Claim \ref{T023}(a) we have $\t_{x''}=\O(\bigcup_{x\in I}\t_x)\subset\t$,
which gives a contradiction.

{\it Case 4}: $\neg \exists\max _{\BL} I$ and $\neg \exists \min _{\BL}F$. Then $(I,F)$ is a gap in $L$
and there is a unique $c\in\widetilde{L}\setminus L$ such that $I<c<F$.
Since $c=\sup I$, by Claim \ref{T023}(a) we have $\t_c=\O(\bigcup_{x\in I}\t_x)\subset\t$
and, since $c=\inf F$, by Claim \ref{T023}(b) we have $\t_c=\bigcap_{x\in F}\t_x\supset\t$.
So $\t_c=\t$ and, hence, $\t _{x_0}\cong \t _c$, where $c\in \widetilde{L}\setminus L$,
which is impossible by Claim \ref{T022}.
\hfill $\Box$
\subsection{The homeomorphism class $[\t_{x_0}]_{\cong}$}
\begin{lem}\label{T036}
If $x\in L$, $\O\in [\t_x]_{\cong}$ and $f:(L_z,\t_x)\rightarrow(L_z,\O)$ is a homeomorphism, then

(a) If $\O\subset\t_x$, then $f\in \Aut (\BL _z)$, $\O=\t_{f(x)}\in \L$ and $f(x)\leq x$;

(b) If $\t_x\subset\O$, then $f\in \Aut (\BL _z)$, $\O=\t_{f(x)}\in \L$ and $f(x)\geq x$.
\end{lem}
\dok
(a) Since $\O\subset\t _x$ and $f[\t _x]=\O$ there are sets $A\subset \widetilde{L}$ and $B\subset(z,x]$ s.t.
\begin{eqnarray}
\O\setminus \{\emptyset,L_z\} & = &\{[z,a)\cap L_z:\ a\in A\}\cup\{(z,b)\cap L_z:\ b\in B\}\label{EQ017}\\
                              & = & \{f[[z,a)\cap L_z]:\ a\in \widetilde{L}\}\cup\{f[(z,b)\cap L_z]:\ b\in(z,x]\}.\nonumber
\end{eqnarray}
First we prove that $f(z)=z$.
By Claim \ref{T020} we have $\{ z\}=\bigcap_{z\in O\in\t_x} O$
and since for $O\subset L_z$ we have $z\in O\in\t_x$ iff $f(z)\in f[O]\in\O$,
we obtain $\{ f(z)\}=f[\{z\}]=f[\bigcap_{z\in O\in\t_x} O]=\bigcap_{f(z)\in f[O]\in\O}\ f[O] \supset\bigcap_{f(z)\in O'\in\O} O'\ni f(z)$;
thus, $\bigcap \{ O'\in\O : f(z)\in O'\}=\{ f(z)\}$.
Assuming that $z<f(z)$,
there would be $e\in (z,f(z))\cap L_z$
and, by the first equality in (\ref{EQ017}),
$e\in \bigcap \{ O'\in\O : f(z)\in O'\}=\{ f(z)\}$, which is false.

Thus $f(z)=z$ and, since $f:L_z \rightarrow L_z$ is a bijection, $f[L]=L$.
For a proof that $f\in \Aut (\BL_z)$ it remains to be shown that
the restriction $f\upharpoonright L :\la L,<\ra\rightarrow\la L,<\ra$ is strictly increasing.
So, if $a_1,a_2\in L$ and $a_1<a_2$,
then $[z,a_2)\cap L_z\in\t_x$ and $a_1\in[z,a_2)\cap L_z\not\ni a_2$;
so we have $f(a_1)\in f[[z,a_2)\cap L_z]\not\ni f(a_2)$.
Since $f[L]=L$ we have $f(a_1),f(a_2)\in L$
and, by (\ref{EQ017}), $f[[z,a_2)\cap L_z]=[z,a)\cap L_z$, for some $a\in A$.
Thus $f(a_1)\in [z,a)\cap L_z\not\ni f(a_2)$,
which implies that $f(a_1)<a\leq f(a_2)$ and $f\in \Aut (\BL_z)$ indeed.

Since $f\in \Aut (\BL_z)$,
like in the proof of  Claim \ref{T018}, there is $F\in\Aut(\widetilde{\BL}_z)$ such that $F\upharpoonright L_z=f$.
Therefore we have $f[[z,a)\cap L_z]=F[[z,a)\cap L_z]=[z,F(a))\cap L_z$, for $a\in \widetilde{L}$,
and $f[(z,b)\cap L_z]=F[(z,b)\cap L_z]=(z,F(b))\cap L_z$, for $b\in (z,x]$.
So, by (\ref{EQ017}),
$\O\setminus \{\emptyset,L_z\}=\{[z,F(a))\cap L_z:\ a\in \widetilde{L}\}\cup\{(z,F(b))\cap L_z:\ b\in(z,x]\}$,
which, since $F(\widetilde{L})=\widetilde{L}$ and $F[(z,x]]=(z,F(x)]=(z,f(x)]$, gives
$\O=\{[z,a'):\ a'\in \widetilde{L}\}\cup\{(z,b'):\ b'\in(z,f(x)]\}\cup\{\emptyset,L_z\}=\t_{f(x)}$.
Since $\O\subset\t_x$ we have $f(x)\leq x$.

(b) Since $\t_x\subset\O$ and $f[\t _x]=\O$ we have
\begin{eqnarray}
\O \setminus \{\emptyset,L_z\} & =       & \{f[[z,a)\cap L_z]:\ a\in \widetilde{L}\}\cup\{f[(z,b)\cap L_z]:\ b\in(z,x]\}\nonumber\\
                               & \supset & \{[z,a)\cap L_z:\ a\in \widetilde{L}\}\cup\{(z,b)\cap L_z:\ b\in(z,x]\}.\label{EQ018}
\end{eqnarray}
First we prove that $f(z)=z$.
Clearly, $z=f(e)$, for some $e\in L_z$,
by (\ref{EQ018}) we have  $\{ z\}=\bigcap_{z\in O\in\O} O$
and since for $O\subset L_z$ we have $z\in O\in\O$ iff $e=f^{-1}(z)\in f^{-1}[O]\in\t _x$,
we obtain $\{ e\}=f^{-1}[\{ z\}]=\bigcap_{e\in f^{-1}[O]\in\t _x} f^{-1}[O]\supset \bigcap_{e\in O'\in\t _x} O'\ni e$;
thus $\bigcap \{ O'\in\t _x : e\in O'\}=\{e\}$
and, by Claim \ref{T020}, $e=z$.

Thus $f(z)=z$, $f[L]=L$ and for a proof that $f\in \Aut (\BL_z)$ we show that
$f\upharpoonright L :\la L,<\ra\rightarrow\la L,<\ra$ is strictly increasing.
On the contrary, suppose that $a_1,a_2\in L$, $a_1<a_2$ and $f(a_1)>f(a_2)$.
Then by (\ref{EQ018}) we have $O:=[z,f(a_1))\cap L_z\in\O$ and, clearly, $f(a_2)\in O\not\ni f(a_1)$;
consequently, $f^{-1}[O]\in\t _x$ and $a_2\in f^{-1}[O]\not\ni a_1$,
which is impossible by (\ref{EQ004}).
So $f\in \Aut (\BL_z)$ and, using (\ref{EQ018}), exactly as in (a) we prove that $\O=\t_{f(x)}$.
Finally, $\t_x\subset\O$ implies that $f(x)\geq x$.
\hfill $\Box$
\begin{prop}\label{T037}
The partial order $\la [\t_{ x_0 }]_{\cong},\subset \ra$ is a union of $2^\k$ disjoint chains isomorphic to $\BL$
such that elements of different chains are incomparable. Thus $\la [\t_{ x_0 }]_{\cong},\subset \ra \cong \dot{\bigcup}_{2^\k}\BL$.
\end{prop}
\dok
Clearly, $[\t_{ x_0 }]_{\cong}=\{ f[\t_{ x_0 }]:f\in \Sym (L_z)\}$.
We recall that $\BL _z =\la L_z ,<\upharpoonright L_z\ra$; for simplicity we will write $<_{L_z}$ instead of $<\upharpoonright L_z$.
Let $f\in \Sym (L_z)$.
If $<_{L_z}^f \subset L_z\times L_z$ is the relation defined by $a<_{L_z}^f b$ iff $f^{-1}(a)<_{L_z}f^{-1}(b)$, then
\begin{equation}\label{EQ024}
\forall a,b\in L_z \;\;(a<_{L_z} b \Leftrightarrow f(a) <_{L_z}^f f(b))
\end{equation}
and, hence, $f:\la L_z ,<_{L_z}\ra \rightarrow \la L_z ,<_{L_z}^f\ra =:\BL _z^f$ is an isomorphism;
thus $\BL _z^f \cong \BL _z$ and $f(z)$ is a smallest element of $\BL _z^f$.

By Proposition \ref{T042} $\L :=\{ \t _x: x\in L \}$ is a maximal chain in $[\t_{ x_0 }]_{\cong}$ isomorphic to $\L$.
Clearly, $\L ^f :=\{ f[\t _x]: x\in L \} \subset [\t_{ x_0 }]_{\cong}$
and if $x_1,x_2\in L$ and $x_1 <x_2$,
then by Claim \ref{T016} $\t _{x_1}\varsubsetneq \t _{x_2}$,
which implies that $f[\t _{x_1}]\varsubsetneq f[\t _{x_2}]$;
so $\la \L ^f, \varsubsetneq\ra\cong \BL$.
If $\L ^f \cup \{ \t\}$ is a chain in $[\t_{ x_0 }]_{\cong}$,
then $\L \cup\{ f^{-1}[\t]\}$ is a chain in $[\t_{ x_0 }]_{\cong}$
and, by the maximality of $\L$, $f^{-1}[\t] =\t _x$, for some $x\in L$;
so, $\t =f[\t _x]\in \L ^f$
and $\L ^f $ is a maximal chain in $[\t_{ x_0 }]_{\cong}$ isomorphic to $\L$.

Clearly we have $\bigcup _{f\in \Sym (L_z)}\L ^f\subset[\t_{ x_0 }]_{\cong}$.
Conversely, if $\t \in [\t_{ x_0 }]_{\cong}$,
then there is a homeomorphism $f:(L_z,\t_{ x_0 })\rightarrow (L_z,\t )$,
which means that $\t =f[\t_{ x_0 }]\in \L ^f$.
So, $[\t_{ x_0 }]_{\cong}=\bigcup _{f\in \Sym (L_z)}\L ^f$ and there are at most $2^{|L|}$-many different chains $\L ^f$.

We show that for $f,g\in \Sym (L_z)$ we have $\L ^f=\L ^g$ or $\L ^f \cap\L ^g=\emptyset$.
So, if $\t\in \L ^f \cap\L ^g$, then $\t =f[\t _{x_1}]=g[\t _{x_2}]$, for some $x_1,x_2\in L$,
and for $x\in L $ we prove that $f[\t _x] \in \L ^g$.
If $x\leq x_1$,
then $\t _x \subset \t _{x_1}$
and, hence, $f[\t _x ]\subset f[\t _{x_1}]=g[\t _{x_2}]$,
which gives $g^{-1}[f[\t _x ]]\subset \t _{x_2}$.
So by Lemma \ref{T036}(a) $g^{-1}[f[\t _x ]]=\t _{x'}$, for some $x'\in L$,
and, hence, $f[\t _x ]=g[\t _{x'}]\in \L ^g$ indeed.
If $x>x_1$, Lemma \ref{T036}(b) gives $f[\t _x] \in \L ^g$ again;
thus $\L ^f \subset \L ^g$ and, by the maximality of $\L ^f$, $\L ^f=\L ^g$.

Assuming that $\L ^f \neq \L ^g $
and that there are comparable $f[\t _{x_1}]\in \L ^f$ and $g[\t _{x_2}]\in \L^g$,
say $f[\t _{x_1}]\subset g[\t _{x_2}]$,
we would have $\t _{x_1}\subset f^{-1}[g[\t _{x_2}]]$
and, by Lemma \ref{T036}(b), $f^{-1}[g[\t _{x_2}]]=\t _x$, for some $x\in L$.
But then $g[\t _{x_2}]=f[\t _x]\in \L ^f \cap\L ^g = \emptyset$, which is false.
Thus, elements of different chains are incomparable.

Further we show that
\begin{equation}\label{EQ022}
\L ^f=\L ^g \Leftrightarrow g^{-1}\circ f \in \Aut (\BL_z) \Leftrightarrow \; <_{L_z}^f=\;<_{L_z}^g.
\end{equation}
If $\L ^f=\L ^g$,
then $f[\t _{x_0}], g[\t _{x_0}] \in \L ^f$
and these topologies are comparable, say $g[\t _{x_0}]\subset f[\t _{x_0}]$,
which gives  $\t _{x_0}\subset g^{-1}[f[\t _{x_0}]]$.
Clearly, the mapping $g^{-1}\circ f: (L_z,\t _{x_0} )\rightarrow (L_z, g^{-1}[f[\t _{x_0}]])$ is a homeomorphism,
by Lemma \ref{T036}(b) we have $g^{-1}\circ f \in \Aut (\la L_z, <_{L_z}\ra)$,
which implies that $g^{-1}[f[<_{L_z}]]=<_{L_z}$,
and, hence, $<_{L_z}^f= f[<_{L_z}]=g[<_{L_z}]=<_{L_z}^g$.

Conversely, if $<_{L_z}^f =<_{L_z}^g$,
that is $f[<_{L_z}]= g[<_{L_z}]$,
then, $g^{-1}[f[<_{L_z}]]= <_{L_z}$
and, hence, $g^{-1}\circ f \in \Aut (\BL_z)$.
Since the automorphism $g^{-1}\circ f$ maps $x_0$ to $g^{-1}(f(x_0))$,
by Claim \ref{T018}(a) the mapping  $g^{-1}\circ f :(L_z,\t _{x_0})\rightarrow(L_z,\t _{g^{-1}(f(x_0))})$ is a homeomorphism,
which gives $g^{-1}[f[\t _{x_0}]]=\t _{g^{-1}(f(x_0))}$
and, hence, $f[\t _{x_0}]=g[\t _{g^{-1}(f(x_0))}]$.
So, since $f[\t _{x_0}]\in \L ^f$ and $g[\t _{g^{-1}(f(x_0))}]\in \L ^g$
we have $\L ^f\cap \L ^g\neq \emptyset$
and, thus,  $\L ^f=\L ^g$.

Finally we show that $|\{ \L ^f:f\in \Sym (L_z)\}|=2^{|L|}$.
By our assumption we have $\k :=|L|=|L_z|\geq \o$; so, $|\Sym (L_z)|= 2^\k$.
Let $\{ P_\xi :\xi <\k\}$ be a partition of the set $L_z$ such that $|P_\xi|=2$, for all $\xi <\k$.
(Identifying $L_z$ with $\k$ we can take $P_\xi=\{ \a _\xi, \a _\xi +1\}$, where $\a _\xi$ is the $\xi$-th even ordinal.)
For a function $\f :\k \rightarrow 2$ let the bijection $f_\f\in \Sym (L_z)$ be defined in the following way.
If $\xi <\k$ and $P_\xi =\{ a_\xi,b_\xi\}$, then
$f_\f (a_\xi)=a_\xi$ and $f_\f (b_\xi)=b_\xi$, if $\f (\xi)=0$;
$f_\f (a_\xi)=b_\xi$ and $f_\f (b_\xi)=a_\xi$, if $\f (\xi)=1$.
By (\ref{EQ022}) it remains to be proved that $\f \neq \p$ implies $<_{L_z}^{f_\f}\neq \;<_{L_z}^{f_\p}$.
So, if $\f \neq \p$ and, say, $\f (\xi) =0$ and $\p (\xi)=1$,
then $f_\f (a_\xi)=a_\xi$, $f_\f (b_\xi)=b_\xi$,
$f_\p (a_\xi)=b_\xi$ and $f_\p (b_\xi)=a_\xi$.
In the initial order $\la L_z,<\ra$ the elements $a_\xi$ and $b_\xi$ are comparable.
If $a_\xi < b_\xi$, then $\la f_\f (a_\xi), f_\f (b_\xi)\ra =\la a_\xi , b_\xi \ra \in \;<_{L_z}^{f_\f}$,
while $\la f_\p (a_\xi), f_\p (b_\xi)\ra =\la b_\xi , a_\xi \ra \in \;<^{f_\p}\setminus <_{L_z}^{f_\f}$.
Thus $<_{L_z}^{f_\f}\neq \;<_{L_z}^{f_\p}$ indeed. If $a_\xi > b_\xi$ we have a similar proof.
\hfill $\Box$
\subsection{The condensation class $[\t_{x_0}]_{\sim}$}
By Proposition \ref{T037} the partial order $\la [\t_{ x_0 }]_{\cong},\subset \ra$
is a union of $2^\k$ disjoint chains of the form $\L ^f:=\{ f[\t _x]: x\in L \}\cong\BL$,
where $f\in S\subset\Sym (L_z)$ and $|S|=2^\k$.
Let
\begin{equation}\label{EQ033}
\Lambda ^f :=\{ f[\t _c]:c\in \widetilde{L}\}, \;\;\mbox{ for }f\in \Sym (L_z).
\end{equation}
By Claim \ref{T016} we have $\la \Lambda ^f , \varsubsetneq\ra \cong \widetilde{\BL}$.
\begin{lem}\label{T040}
$\Conv (\L ^f)=\Lambda^f$, for each $f\in \Sym (L_z)$.
\end{lem}
\dok
We recall that $\L :=\{ \t _x: x\in L \}$, $\Lambda :=\{ \t _c:c\in \widetilde{L}\}$
and first prove that $\Conv (\L)=\Lambda$.
If $\O \in \Conv (\L)$,
then $\t _{x_1}\subset \O \subset \t _{x_2}$, for some $x_1,x_2\in L$,
and, hence, $\O \setminus \t _{x_1}=\{(z,b)\cap L_z:\ b\in B\}$, where $B\subset (x_1,x_2]$.
If $B=\emptyset$, then $\O =\t _{x_1}\in \L \subset \Lambda$.
Otherwise, there is $c=\sup_{\widetilde{L}}B \leq x_2$.
If $c\in B$, we have $(z,c)\cap L_z\in\O$;
and if $c\not\in B$, then for each $x\in (z,c)\cap L_z$ there is $b\in B$ such that $x\in (z,b)\cap L_z\in \O$
and, hence, $(z,c)\cap L_z\in\O$ again.
Since $\t_{\leftarrow}\subset\O$ for each $b\leq c$ we have $(z,c)\cap L_z\cap [z,b)\cap L_z=(z,b)\cap L_z\in \O$;
so, $B=(x_1,c]$ and, hence, $\O=\t_c\in\Lambda$.
Conversely, if $\t_c\in\Lambda$ (where $c\in \widetilde{L}$)
then, since $L$ is dense in $\widetilde{L}$ and has no end points,
there are $x_1,x_2\in L$ such that $x_1<c <x_2$.
By Claim \ref{T016} we have $\t_{x_1}\subsetneq\t_c\subsetneq\t_{x_2}$
and, hence, $\t _c \in \Conv (\L)$.

Now, if $\O \in\Conv (\L ^f)$,
then $f[\t _{x_1}]\subset \O \subset f[\t _{x_2}]$, for some $x_1,x_2\in L$,
and, hence,  $\t _{x_1}\subset f^{-1}[\O] \subset \t _{x_2}$,
which implies that $f^{-1}[\O]\in \Conv (\L)=\Lambda$.
So, $f^{-1}[\O]=\t _c$, for some $c\in \widetilde{L}$
and $\O=f[\t _c]\in \Lambda^f$.
Conversely, if $f[\t _c]\in \Lambda^f$,
then $\t _c\in \Lambda =\Conv (\L)$
and, hence, $\t _{x_1}\subset \t_c \subset \t _{x_2}$, for some $x_1,x_2\in L$,
Consequently $f[\t _{x_1}]\subset f[\t_c ] \subset f[\t _{x_2}]$
and, since $f[\t _{x_1}], f[\t _{x_2}]\in\L ^f$,
we have $f[\t _c]\in\Conv (\L ^f)$.
\hfill $\Box$
\begin{prop}\label{T041}
The partial order $\la [\t_{ x_0 }]_{\sim},\subset \ra$ is a union of $2^\k$ disjoint chains isomorphic to $\widetilde{\BL}$
such that elements of different chains are incomparable. Thus $\la [\t_{ x_0 }]_{\sim },\subset \ra \cong \dot{\bigcup}_{2^\k}\widetilde{\BL}$.
\end{prop}
\dok
By Proposition \ref{T037} we can choose bijections $f_\a \in \Sym (L_z)$, for $\a < 2^\k$, such that $[\t _{x_0}]_{\cong}=\dot{\bigcup}_{\a <2^\k}\L ^{f_\a}$.
Since $[\t _{x_0}]_{\sim}= \Conv ([\t _{x_0}]_{\cong})$ first we show that
\begin{equation}\label{EQ032}\textstyle
\Conv(\bigcup_{\a <2^\k}\L ^{f_\a})=\bigcup_{\a <2^\k}\Lambda ^{f_\a}.
\end{equation}
If $\O \in \Conv(\bigcup_{\a <2^\k}\L ^{f_\a})$,
then $\s \subset \O \subset \t$, for some $\s ,\t \in \bigcup_{\a <2^\k}\L ^{f_\a}$,
and, since by Proposition \ref{T037} the elements of different chains $\L ^{f_\a}$ are incomparable,
there is $\a \in 2^\k$ such that $\s,\t\in \L ^{f_\a}$.
So, by Lemma \ref{T040} we have $\O \in \Conv (\L ^{f_\a})=\Lambda^{f_\a}$.
Conversely, for $\a<2^\k$ we have $\L ^{f_\a}\subset \bigcup_{\a <2^\k}\L ^{f_\a}$;
so, by Lemma \ref{T040}, $\Lambda^{f_\a} =\Conv (\L ^{f_\a})\subset \Conv(\bigcup_{\a <2^\k}\L ^{f_\a})$
and (\ref{EQ032}) is true.

Second, for different $\a ,\b<2^\k$, $\O \in \Lambda^{f_\a}$ and $\O' \in \Lambda^{f_\b}$,
we show that the topologies $\O$ and $\O '$ are incomparable.
Assuming that $\O \subset \O '$
there would be $\t \in \L^{f_\a}$ and $\s \in \L^{f_\b}$ such that $\t \subset \O$ and $\O '\subset \s$
and, hence, $\t \subset \s$, which is impossible because by Proposition \ref{T037} the elements of $\L^{f_\a}$ and $\L^{f_\b}$ are incomparable.
Consequently we have $\Lambda^{f_\a}\cap \Lambda^{f_\b}=\emptyset$
and, hence, $[\t _{x_0}]_{\sim}=\dot{\bigcup}_{\a <2^\k}\Lambda ^{f_\a}$.
Since $\Lambda^{f_\a}\cong \widetilde{\BL}$, for all $\a <2^\k$, we are done.
\hfill $\Box$
\subsection{The remainder $[\t_{ x_0 }]_{\sim}\setminus [\t_{ x_0 }]_{\cong}$}
By Claim \ref{T018}(c) and Lemma \ref{T040}
$\Lambda =\{ \t _c :c\in \widetilde{L}\}=\Conv (\L)=\widetilde{\L}\cong \widetilde{\BL}$ is a maximal chain in $[\t_{ x_0 }]_{\sim}$
containing $\t_{ x_0 }$.
Let $\thickapprox$ be the equivalence relation on $\widetilde{L}$
defined by $c\thickapprox c'$ iff $\t _c \cong \t _{c'}$
and let $(\widetilde{L}/\!\thickapprox) \setminus  \{ [x_0 ]_{\thickapprox}\} =\{[c_\xi]_{\thickapprox }:\xi <\l\}$
be an enumeration of its classes different from $[x_0 ]_{\thickapprox}$.
Then, for $\xi <\l$, defining $L_\xi :=[c_\xi]_{\thickapprox}$,
we obtain a suborder $\BL _{\xi}:= \la L_\xi, <\ra$ of $\widetilde{\BL}$
and we have $\widetilde{L}=L\cup \bigcup _{\xi <\l}L_\xi$.
\begin{prop}\label{T048}
For  $\xi <\l$ we have $[\t_{ c_\xi }]_{\sim}=[\t_{ x_0 }]_{\sim}$
and $\la [\t _{c_\xi}]_{\cong}, \subset\ra \cong \dot{\bigcup}_{\a <2^\k}\BL _{\xi}$.
\end{prop}
\dok
Clearly, $\L _\xi :=\{ \t _c :c\in L_\xi\}=\Lambda \cap [\t _{c_\xi}]_{\cong }$
and $\Lambda =\L \cup \bigcup _{\xi <\l}\L _\xi$.
By Claim \ref{T016}, for $f\in \Sym (L_z)$ we have $\Lambda ^f :=\{ f[\t _c]:c\in \widetilde{L}\}\cong \widetilde{\BL}$
and, defining $\L _\xi ^f :=\{ f[\t _c] :c\in L_\xi\}$, we have $\Lambda ^f =\L ^f \cup \bigcup _{\xi <\l}\L _\xi ^f$.
So, by Proposition \ref{T041} and (\ref{EQ032})
\begin{equation}\label{EQ035}\textstyle
[\t _{x_0}]_{\sim }=\bigcup_{\a <2^\k}\Lambda ^{f_\a}=\bigcup_{\a <2^\k}(\L ^{f_\a} \cup \bigcup _{\xi <\l}\L _\xi ^{f_\a}).
\end{equation}
If $\xi <\l$, then by (m1) there are $x,x'\in L$ such that $x<c_\xi<x'$,
which by Claim \ref{T016} gives $\t _x \varsubsetneq \t _{c_\xi} \varsubsetneq \t _{x'}$,
and, hence, $\t _x \preccurlyeq \t _{c_\xi} \preccurlyeq \t _{x'}$ and $\t _{c_\xi} \sim \t _{x'}$.
Thus $[\t_{ c_\xi }]_{\sim}=[\t_{ x_0 }]_{\sim}$ and $[\t _{c_\xi}]_{\cong} \subset [\t_{ x_0 }]_{\sim}$,
which by (\ref{EQ035}) gives $[\t _{c_\xi}]_{\cong}=\bigcup_{\a <2^\k}\L _\xi ^{f_\a}$.
Since $\la \L _\xi ^{f_\a},\varsubsetneq\ra\cong\la \L _\xi ,\varsubsetneq\ra\cong$
and by Claim \ref{T016} we have $\la \L _\xi ,\varsubsetneq\ra\cong \BL _{\xi}$,
we finally obtain $\la [\t _{c_\xi}]_{\cong}, \subset\ra \cong \dot{\bigcup}_{\a <2^\k}\BL _{\xi}$.
\kdok
Thus, if $\widetilde{\BL}\neq \BL$, then $|\,[\t _{x_0}]_{\sim}/\!\cong|=\l +1\geq 2$
and, coloring homeomorphic topologies in the same color,
the condensation class $[\t _{x_0}]_{\sim}$ is a union of $2^\k$-many maximal chains $\Lambda ^{f_\a}\cong \widetilde{\BL}$, $\a <2^\k$,
colored in $(\l+1)$-many colors and isomorphic as colored linear orders.
We describe one situation with 2 colors.
\begin{prop}\label{T047}
If $\BL$ and $\widetilde{\BL}\setminus \BL$
are 1-homogeneous linear orders
and each $f\in \Aut (\widetilde{\BL}\setminus \BL)$ extends to an automorphism of  $\widetilde{\BL}$,
\footnote{in particular, if $\widetilde{\BL}$ is a Dedekind completion of $\widetilde{\BL}\setminus \BL$ too; then by (m2) $\BL$ and $\widetilde{\BL}\setminus \BL$ are dense.}
then taking arbitrary $c\in \widetilde{L}\setminus L$ we have
$[\t _{x_0}]_{\sim}=[\t _{x_0}]_{\cong }\cup [\t _c]_{\cong}\cong \dot{\bigcup} _{2^\k}\widetilde{\BL}$
and $\la [\t _c]_{\cong}, \subset \ra \cong \dot{\bigcup} _{2^\k}(\widetilde{\BL}\setminus \BL)$.
\end{prop}
\dok
By Proposition \ref{T048} it remains to be proved that
the topologies $\t _c$, $c\in \widetilde{L}\setminus L$, are homeomorphic.
We note that we can not simply replace $\t _{x_0}$ by $\t _c$ in Claim \ref{T018}(a),
because $\t _c $ is a topology on the set $\{ z\} \cup L$, not on the set $\{ z\} \cup (\widetilde{L}\setminus L)$.

By the assumptions, for $c_0,c_1\in \widetilde{L}\setminus L$
there are $\f\in \Aut (\widetilde{\BL}\setminus \BL)$ such that $\varphi(c_0)=c_1$
and $F\in\Aut(\widetilde{\BL}_z )$ such that $F|(\widetilde{L}\setminus L)=\f$ and $F(z)=z$.
Now, $f:=F|L_z \in \Aut (\BL _z)$
and (see the proof of Claim \ref{T018}(a))
$f:(L_z ,\t _x)\rightarrow (L_z ,\t _{f(x)})$ is a homeomorphism;
that is, $f[\t _x]=\t _{f(x)}$, for all $x\in L$.
For $i<2$ we have $c_i=\inf _{\widetilde{\BL}}((c_i,\cdot)\cap L)$
and, by Claim \ref{T023}(b), $\t _{c_i}=\bigcap _{x\in (c_i,\cdot)\cap L}\t _x$.
Thus, since $f[(c_0,\cdot)\cap L]=F[(c_0,\cdot)\cap L]=(c_1,\cdot)\cap L$,
we have $f[\t _{c_0}]=\bigcap _{x\in (c_0,\cdot)\cap L}f[\t _x]=\bigcap _{x\in (c_0,\cdot)\cap L}\t _{f(x)}=\bigcap _{y\in (c_1,\cdot)\cap L}\t _y=\t _{c_1}$.
\hfill $\Box$
\begin{ex}\label{EX001}\rm
{\it Applications of Proposition \ref{T047}.}
Taking $\BL =\Z\Z$ it is easy to check that $\widetilde{\BL}\cong (\Z +1)\Z$, $\widetilde{\BL}\setminus \BL\cong \Z$
and each automorphism of $\widetilde{\BL}\setminus \BL$ extends to an automorphism of  $\widetilde{\BL}$.
Thus taking $x_0\in L$ and $c\in \widetilde{L}\setminus L$
we have $[\t _{x_0}]_{\sim}=[\t _{x_0}]_{\cong }\cup [\t _c]_{\cong}\cong \dot{\bigcup} _{\fc}(\Z +1)\Z$,
where $[\t _{x_0}]_{\cong} \cong \dot{\bigcup} _{\fc}\Z\Z$
and $\la [\t _c]_{\cong}, \subset \ra \cong \dot{\bigcup} _{\fc}\Z$.

Taking $\BL =\Q$ we obtain topologies $\t _0$ and $\t _{\pi}$ on a {\it countable} set
such that $[\t _0]_{\sim}=[\t _\pi]_{\sim}\cong \dot{\bigcup}_{\fc}\R$, $[\t _0]_{\cong}\cong \dot{\bigcup}_{\fc}\Q$ and $[\t _\pi ]_{\cong}\cong \dot{\bigcup}_{\fc}(\R\setminus\Q)$.
On the other hand, taking $\BL =\R\setminus\Q$ we obtain topologies $\t _0'$ and $\t _{\pi}'$ on a {\it set of size $\fc$}
such that $[\t _0']_{\sim}\cong \dot{\bigcup}_{2^\fc}\R$, $[\t _0']_{\cong}\cong \dot{\bigcup}_{2^\fc}\Q$ and $[\t _\pi ']_{\cong}\cong \dot{\bigcup}_{2^\fc}(\R\setminus\Q)$.

Moreover, under CH there are $\fc ^+$-many non-isomorphic new  examples.
First, Adams \cite{Ada} proved that, under CH, whenever $C_\a \in [\R]^{\fc}$, for $\a <\fc$,
there is a $2$-homogeneous linear order $\BL\subset \R$ such that $\BC_\a \not\hookrightarrow \BL$, for all $\a <\fc$.
Second, if $\BL \subset \R$ is $2$-homogeneous and dense,
then $\R \setminus \BL $ is $2$-homogeneous (see \cite{Ros}, p.\ 163);
so, since Adams' order $\BL$ is dense (it contains $\Q$, see \cite{Ada}, p.\ 160)
$\R \setminus \BL$ is $2$-homogeneous.
Third, requesting that e.g.\ $(0,1)\not\hookrightarrow \BL$ we obtain that $\R \setminus \BL$ is dense in $\R$,
and, hence, $\R$ is a Dedekind completion of $\R \setminus \BL$.
Thus the linear order $\BL$ satisfies the assumptions of Proposition \ref{T047}.
Now, an easy recursion gives non-isomorphic linear orders $\BL _\a$, $\a <\fc ^+$, satisfying these assumptions.
So, for $\a <\fc ^+$, taking $x_0\in L_\a$ and $c\in \R \setminus L_\a$ we have
$[\t _{x_0}]_{\sim}=[\t _c]_{\sim}\cong \dot{\bigcup}_{2^\fc}\R$,
$[\t _{x_0}]_{\cong}\cong \dot{\bigcup}_{2^\fc}\BL _\a$ and
$[\t _c ]_{\cong}\cong \dot{\bigcup}_{2^\fc}(\R\setminus\BL _\a)$.
\end{ex}
\section{Weakly reversible, non-reversible topologies}\label{S4}
If $\BL $ is a $1$-homogeneous linear order which is not Dedekind complete,
then by Claim \ref{T022} we have $[\t _{x_0} ]_{\cong} \varsubsetneq [\t_{x_0} ]_{\sim}$;
so, the topology $\t_{x_0}$ is not weakly reversible.
If $\BL $ is Dedekind complete,
we have the following consequence of  Theorem \ref{T003}.
\begin{te}\label{T043}
For each $1$-homogeneous Dedekind complete linear order $\BL $ of size $\k\geq \o$ there is a weakly reversible non-reversible topology $\t \in \Top _\k$ such that
\begin{equation}\label{EQ036}\textstyle
\la [\t ]_{\cong}, \subset \ra =\la [\t ]_{\sim}, \subset \ra \cong \dot{\bigcup}_{2^\k}\BL
\end{equation}
and, hence, each maximal chain in $[\t ]_{\cong}$ is isomorphic to $\BL$.
\end{te}
Concerning the application of Theorem \ref{T043}
we note that $\Z$ is, up to isomorphism,
the unique scattered $1$-homogeneous Dedekind complete linear order,
while the non-scattered ones are dense (see \cite{Ohk}) and $\R$ is an archetypical example.
Using the following elementary statement we obtain more dense examples.
\begin{fac}\label{T045}
If $\BL$ is a dense $1$-homogeneous Dedekind complete linear order
such that $(a,b)\cong \BL$, whenever $a<b\in L$,
then $\BL$ is $2$-homogeneous and the following three linear orders are $2$-homogeneous and Dedekind complete
\begin{eqnarray*}
R(\BL) & := & \BL + (1+\BL)\o _1, \\
L(\BL) & := & (\BL +1)\o _1^* +\BL , \\
M(\BL) & := & (\BL +1)\o _1^* +\BL+ (1+\BL)\o _1.
\end{eqnarray*}
In addition, $\cf (\BL)=\ci (\BL)=\chi (x)=\o$, for all $x\in L$.
\end{fac}
\dok
Assuming that $\BL$ contains a bounded copy of $\o _1$,
by completeness there would be points of character $\o$ and $\o _1$,
which is impossible by homogeneity.
Thus $\chi (x)=\o$, for all $x\in L$,
and, since $(a,b)\cong \BL$, we have $\cf (\BL)=\ci (\BL)=\o$.

Taking $\o$-sequences $a=x_0 <x_1 < \dots$, cofinal in $(a,b)$,
and $a=y_0 <y_1 < \dots$, cofinal in $(a,\cdot)$,
we have $(a,b)=\bigcup _{n\in \o}(x_n,x_{n+1}]$ and $(a,\cdot )=\bigcup _{n\in \o}(y_n,y_{n+1}]$;
so, since $(s,t]\cong \BL +1$, whenever $s<t\in L$,
we have $(a,\cdot)\cong \BL$ and, similarly, $(\cdot ,a)\cong \BL$, for all $a\in L$.
Now, if $a<b$ and $c<d$,
it is easy to make $f\in \Aut (\BL)$ such that $f(a)=c$ and $f(b)=d$;
thus $\BL$ is $2$-homogeneous indeed.

We show that $R(\BL)$ is Dedekind complete.
First, $R(\BL)= \sum _{\a <\o _1}\BL _\a$, where $\BL _0 \cong \BL$ and $\BL _\a \cong 1+ \BL$, for $\a >0$.
If $\emptyset \neq X\subset R(\BL)$ has an upper bound,
let $\b :=\min \{ \a <\o _1: L _\a \mbox{ contains an upper bound for }X \}$.
Then the set $X \cap \bigcup _{\a <\b}L _\a$ is unbounded in $\bigcup _{\a <\b}L _\a$.
If $y:= \min \BL _\b $ and $X\cap L_\b \subset \{y\}$, then $\sup X =y$;
otherwise, the set $X\cap L_\b$ has an upper bound in $\BL _\b \setminus \{y\}\cong \BL$;
since $\BL$ is Dedekind complete there is $z:=\sup _{\BL _\b \setminus \{y \}}(X\cap L_\b)$
and, clearly, $z=\sup _{R(\BL)}X$.

We prove that $\Lambda _\b :=\sum _{\a <\b}\BL _\a\cong \BL$, for each limit ordinal $\b <\o _1$.
Since $\BL$ is dense it contains a bounded copy of $\b$,
say $\{ x_\a :\a <\b\}$.
Then, since $(x_\a,x_{\a +1})\cong \BL$, for $\a <\b$,
we have $X:=(x_0,x_1)\cup\bigcup _{0<\a <\b}[x_\a,x_{\a +1})\cong \Lambda _\b$.
For $x:=\sup \{ x_\a :\a <\b\}$ we have $X=(x_0 ,x)\cong \BL$
and, hence, $\Lambda _\b\cong \BL$.

In $R(\BL) $, let $x_0<x_1$ and $y_0<y_1$
and let $\b <\o _1$ be a limit ordinal such that $x_0,x_1,y_0,y_1\in \Lambda _\b$
Since $\Lambda _\b\cong \BL$ there is $f\in \Aut (\Lambda _\b)$ such that $f(x_i)=y_i$, for $i<2$,
and, clearly, $f\cup \id _{\sum _{\b \leq \a <\o_1}\BL _\a}\in \Aut (R(\BL))$;
thus, $R(\BL) $ is a $2$-homogeneous linear order.
For $L(\BL)=R(\BL ^*)^*$ we have a dual proof. For $M(\BL)$ we proceed similarly;
here we have $(\BL +1)\b ^* +\BL + (1+\BL )\b \cong\BL$.
\kdok
By the following result of Ohkuma \cite{Ohk}, each dense $1$-homogeneous Dedekind complete linear order
is generated by one satisfying the assumptions of Fact \ref{T045}.
\begin{fac}{(\cite{Ohk})}\label{T046}
If $\BL$ is a dense $1$-homogeneous Dedekind complete linear order,
then there is a dense $1$-homogeneous Dedekind complete linear order $\BC$
such that $(a,b)\cong \BC$, whenever $a<b\in C$,
and $\BL$ is isomorphic to $\BC$, $R(\BC)$, $L(\BC)$ or $M(\BC)$.
\end{fac}
We note that for linear orders $\BL$ and $\BC$ from Fact \ref{T046} we have $|L|=|C|=\fc$.
Namely, in $\BL$, a closed interval $[a,b]$, where $a<b$, is a complete linear order;
so, it is a compact space (\cite{Eng}, p.\ 221).
Since, by Fact \ref{T045}, $\chi (x)=\o$, for all $x\in [a,b]$,
we have $|[a,b]|= \fc$ (\cite{Hod}, p.\ 30, 31).
Thus, since $\cf (\BL),\ci (\BL)\leq \o _1$,
the set $L$ is a union of $\leq \o _1$ intervals of size $\fc$
and, hence, $|L|=|C|=\fc$.
\begin{ex}\label{R000}\rm
{\it Applications of Theorem \ref{T043}.}
By Fact \ref{T045}, the linear orders $\BC$ and $\BL$ from Fact \ref{T046} are, in fact, $2$-homogeneous;
thus, $(x,y)\cong \BC$, for $x<y\in L$.
So, if some of these intervals is separable, then $\BC \cong \R$
and, taking  $\BL$ to be $\R$ or $R(\R)$ (the ``long line" without the first element),
or $L(\R)$ (its reverse) or $M(\R)$,
by Fact \ref{T045} and Theorem \ref{T043} we obtain a topology $\t \in \Top _\fc$ satisfying (\ref{EQ036}), for $\k =\fc$.

In addition, there are such $\BC$-s with no separable intervals.
For example, Shelah in \cite{She} constructed (in ZFC) a Countryman line
and sketched a construction of a 2-homogeneous {\it Aronszajn continuum}
(a complete dense non-separable first countable linear order $\BA$
such that $\o _1,\o _1^* \not \hookrightarrow \BA$ and $w(\overline{X})=\o$, for each $X\in [A]^\o$);
for a construction see the chapter of Todor\v cevi\'c, \cite{Tod}, p.\ 260.
Thus, taking for $\BC$ the order $\BA$ without end points,
since $\cf(\BC)=\ci(\BC)=\chi (x)=\o$, for all $x\in C$,
and all the intervals $(a,b)$ in $\BC$ are isomorphic,
partitioning $(a,b)$ and $C$ into intervals
as in the proof of Fact \ref{T045} we show that $(a,b)\cong \BC$, whenever $a<b\in C$.
By the same fact $R(\BC)$, $L(\BC)$ and $M(\BC)$ are $2$-homogeneous and Dedekind complete linear orders
and, clearly, have no separable intervals.
So, for $\BL\in \{\BC ,R(\BC), L(\BC), M(\BC)\}$ by Theorem \ref{T043} there is $\t \in \Top _\fc$ satisfying (\ref{EQ036}).

Moreover in \cite{Mil} Hart and van Mill constructed $2^{\fc}$ non-isomorphic $2$-homog.\ continua
such that each uncountable suborder contains an uncountable real type;
thus these continua are not Aronszajn.
If $\BC$ is a linear order obtained from some of these continua by deleting of its end points,
then $\cf (\BC)=\ci (\BC)=\o$
and  for $\BL\in \{\BC ,R(\BC), L(\BC), M(\BC)\}$ by Theorem \ref{T043} there is $\t \in \Top _\fc$ satisfying (\ref{EQ036}).
\end{ex}
\paragraph{Acknowledgement}
The research of both authors was supported by the Science Fund of the Republic of Serbia,
Program IDEAS, Grant No.\ 7750027:
{\it Set-theoretic, model-theoretic and Ramsey-theoretic
phenomena in mathematical structures: similarity and diversity}--SMART.

\end{document}